\documentclass[12pt]{amsart}
\usepackage{amsfonts}
\usepackage{amssymb,mathrsfs}
\usepackage{amsmath}
\usepackage{enumerate}
\usepackage{latexsym}
\usepackage{graphicx}
\usepackage{bm}
\usepackage{subfigure}
\usepackage{float}
\usepackage{amsmath}
\usepackage{amsthm}
\usepackage{verbatim}
\usepackage{vmargin}
\usepackage{amstext}
\usepackage{array}
\usepackage{booktabs}
\usepackage{indentfirst}
\usepackage[compress]{cite}
\usepackage{url}
\usepackage{geometry}
\numberwithin{equation}{section}
\newtheorem{thm}{Theorem}[section]
\newtheorem{cor}[thm]{Corollary}
\newtheorem{lma}[thm]{Lemma}

\newtheorem{defn}[thm]{Definition}
\newtheorem{rem}[thm]{Remark}

\theoremstyle{definition}

\begin{document}
\setlength{\oddsidemargin}{80pt}
\setlength{\evensidemargin}{80pt}
\setlength{\topmargin}{80pt}
\author{Jun~Xian, Xiaoda~Xu}

\address{J.~Xian\\Department of Mathematics and Guangdong Province Key Laboratory of Computational Science
 \\Sun
Yat-sen University
 \\
510275 Guangzhou\\
China.} \email{xianjun@mail.sysu.edu.cn}

\address{X.~Xu\\Department of Mathematics
 \\Sun
Yat-sen University
 \\
510275 Guangzhou\\
China.} \email{xuxd26@mail2.sysu.edu.cn}

\title[Expected $L_2$-discrepancy for stratified sampling]{EXPECTED $L_2-$DISCREPANCY BOUND FOR A CLASS OF NEW STRATIFIED SAMPLING MODELS}


\date{\today}

\begin{abstract}
We introduce a class of convex equivolume partitions. Expected $L_2-$discrepancy are discussed under these partitions. There are two main results. First, under this kind of partitions, we generate random point sets with smaller expected $L_2-$discrepancy than classical jittered sampling for the same sampling number. Second, an explicit expected $L_2-$discrepancy upper bound under this kind of partitions is also given. Further, among these new partitions, there is optimal expected $L_2-$discrepancy upper bound.
\end{abstract}

\maketitle
\section{Introduction}\label{intro}

In real sampling processes, it is necessary to know how well-spread these sampling points are. One can select the sampling set randomly which has achieved successful applications in the field of Monte Carlo simulation, compressed sensing, image processing and learning theory \cite{CanTao,chanzick,cucksmale,FrenKJS,liangf,wwhlf,Kriegvca}. The concept of discrepancy is a fundamental building block in the quantification of many point distributions problems. There is a list of interesting discrepancy measures, such as star discrepancy, extreme discrepancy, $G-$discrepancy, isotrope discrepancy, lattice discrepancy, and so on (see e.g., \cite{DGW2014, HEFP}). Among them, $L_{2}-$discrepancy is the most widely studied.

\textbf{$L_{2}-$discrepancy}. $L_{2}-$discrepancy of a sampling set $P_{N, d}=\{t_{1}, t_{2}, \ldots , t_{N}\}$ is defined by

\begin{equation}\label{lpdefn}
L_{2}(D_{N},P_{N, d})=\Big(\int_{[0,1]^{d}}|\lambda([0,z))- \frac{1}{N}\sum_{i=1}^{N}\mathbf{1}_{[0,z)}(t_{i})|^{2}dz\Big)^{1/2},
\end{equation}
where $\lambda$ denotes the Lebesgue measure, $\mathbf{1}_{A}$ denotes the characteristic function on set $A$. For the applications of $L_{2}-$discrepancy, see\cite{Dick2005,Dick2006,Dick2014,Dick2020}.

In the definition of $L_2-$discrepancy, if we introduce the counting measure $\#$, \eqref{lpdefn} can also be expressed as 

\begin{equation}
    L_{2}(D_{N},P_{N, d})=\Big(\int_{[0,1]^{d}}|\lambda([0,z))-
\frac{1}{N}\#\big(P_{N, d}\cap[0,z)\big)|^{2}dz\Big)^{1/2},
\end{equation}
where $\#\big(P_{N, d}\cap[0,z)\big)$ denotes the number of points falling into the set $[0,z).$

To simplify the expression of $L_{2}-$discrepancy, we employ the discrepancy function $\Delta(P_{N,d},z)$ via:

\begin{equation}\label{disfun1}
    \Delta(P_{N,d},z)=\lambda([0,z))-
\frac{1}{N}\#\big(P_{N, d}\cap[0,z)\big).
\end{equation}

Accordingly, the $L_2-$discrepancy can be extended to a fixed compact convex set $K\subset \mathbb{R}^{d}$ with $\lambda(K)>0,$ see \cite{KP}. Discrepancy function in \eqref{disfun1} of a finite set of points $P=\{x_1,x_2,\ldots,x_n\}\subset K$ is now given by

\begin{equation}\label{deltapx}
    \Delta(P,x)=\frac{\lambda\big((-\infty,x]\cap K\big)}{\lambda(K)}-
\frac{1}{N}\#\big(P\cap(-\infty,x]\big).
\end{equation}

For fixed $d$, the best known asymptotic upper bounds for discrepancy are of the form $$O(\frac{(\ln N)^{\alpha_{d}}}{N}),$$ where $\alpha_{d}\ge 0$ are constants depending on dimension $d$. These involve special deterministic point set constructions, which are \textbf{low discrepancy point sets}. Examples of such point sets can be found in \cite{HNie,Dick2010QMC}. For applications arising in computer graphics, quantitative finance and learning theory, see e.g., \cite{Cristiano,Blue-Noise,rankla,phdpaper}.

Although low discrepancy (deterministic) point sets are widely used in numerical integration, the simulation of many phenomena in the real world requires the introduction of random factors. Recently, a large amount of research investigating random sampling for different function spaces has emerged in \cite{BG1,BG2,FX2019}, due to the simplicity, flexibility and effectiveness of the subject. Besides, in the field of discrepancy, probabilistic star discrepancy bounds for Monte Carlo point sets are considered in \cite{aistleitner2014,Gnewuch2020}, while centered discrepancy of random sampling and Latin hypercube random sampling are investigated in \cite{FMW2002}. Motivated by these developments, we incorporate a random viewpoint into our study of star discrepancy to consider a special random sampling method, which is \textbf{stratified sampling}. Its special case is called \textbf{jittered sampling} that is formed by grid-based equivolume partition.

Some random sampling strategies, for example, simple random sampling, stratified sampling, Latin hypercube sampling, etc. are commonly used in the real sampling process, see \cite{shirely1994,McBec,stein87}. Formers have made sufficient research on estimating the expected discrepancy with random samples. For researches on expected star discrepancy of jittered sampling, we refer to \cite{jittsamp,Doerr2}. Both the upper and the lower bounds for the discrepancy of jittered sampling are given in \cite{jittsamp}, while the bounds in \cite{Doerr2} improve them and remove the asymptotic requirement that $m$ is sufficiently large compared to dimensions $d$(where $N=m^d$ means the number of subcubes of grid-based equivolume partition). Starting from the discrepancy itself, rather than estimating its bound. In \cite{KP2}, it is shown that jittered sampling construction gives rise to a set whose expected $L_p-$discrepancy is smaller than that of purely random points. Further, a theoretical conclusion that the jittered sampling does not have the minimal expected $L_2-$discrepancy among all stratified samples from convex equivolume partitions with the same number of points is presented in \cite{KP}. Our research will be carried out on the $d$-dimensional unit cube, which can be easily extended to a more general compact convex set. Studies on convex bodies are extensive, see \cite{BGK,GKM}. In the following, we shall construct a class of convex body partitions to analyze expected $L_2-$discrepancy, which turns out to provide better results than jittered sampling.

Throughout this paper, we adopt the idea of stratified random sampling to study $L_2-$discrepancy. First, we design an infinite family of partitions with partition parameter $0\leq\theta\leq\frac{\pi}{2}$ that generates point sets with a smaller expected $L_2-$discrepancy than classical stratified sampling for sampling number $N=m^d$, which is, $$\mathbb{E}(L_2^2(D_N,P_{\Omega^{*}_{\sim}}))\leq\mathbb{E}(L_2^2(D_N,P_{\Omega^{*}_{|}})),$$where $P_{\Omega^{*}_{\sim}}$ and $P_{\Omega^{*}_{|}}$ denote stratified samples generated by the new infinite family of partitions and grid-based equivolume partition respectively. The equal signs hold if and only if stratified sampling sets $P_{\Omega^{*}_{\sim}}$ are selected for jittered sampling set $P_{\Omega^{*}_{|}}$. 
Second, optimal expected $L_2-$discrepancy bound is also provided under this class of partitions. That is, they are better than the employment of jittered sampling. We obtain the following explicit estimation $$\mathbb{E}(L_2^2(D_N,P_{\Omega^{*}_{\sim}}))\leq \frac{d}{N^{1+\frac{1}{d}}}+\frac{1}{N^3}\cdot\frac{1}{3^{d-2}}\cdot P(\theta),$$ where $P(\theta)$ is the function about partition $\theta$. Taking $\theta=arctan\frac{1}{2}$ and $\theta=0$, we can obtain the upper bounds for optimal partition and grid-based equivolume partition respectively. 

The rest of this paper is organized as follows. In Section \ref{prelim} we present some preliminaries on stratified sampling and newly designed partition models. In Section \ref{sec3} we provide comparisons of the expected $L_2-$discrepancy for stratified sampling under a kind of convex equivolume partitions. The explicit expected $L_2-$discrepancy upper bounds for these newly stratified models are also obtained. In Section \ref{pfmr} we include the proofs of all theorems and lemmas. Finally, in Section \ref{conclu} we conclude the paper. 

\section{Preliminaries on stratified sampling and new partition models}\label{prelim}

Before introducing the main result, we list preliminaries used in this paper. 

\subsection{Stratified sampling}
\begin{figure*}[h]
\centering
\subfigure[two dimensional case]{
\begin{minipage}{7cm}
\centering
\includegraphics[width=0.6\textwidth]{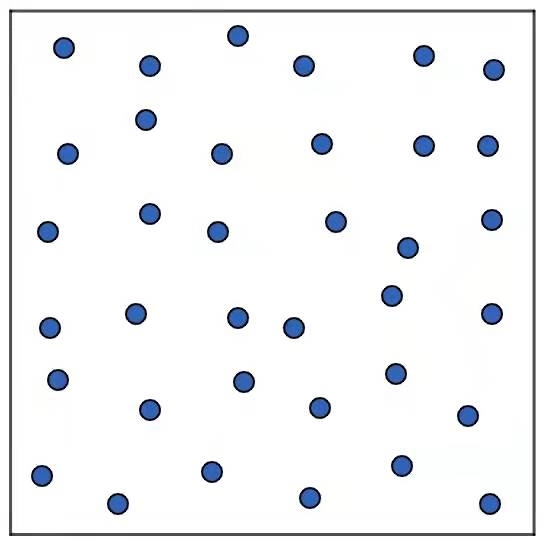}
\end{minipage}
}
\subfigure[three dimensional case]{
\begin{minipage}{5cm}
\centering
\includegraphics[width=1.0\textwidth]{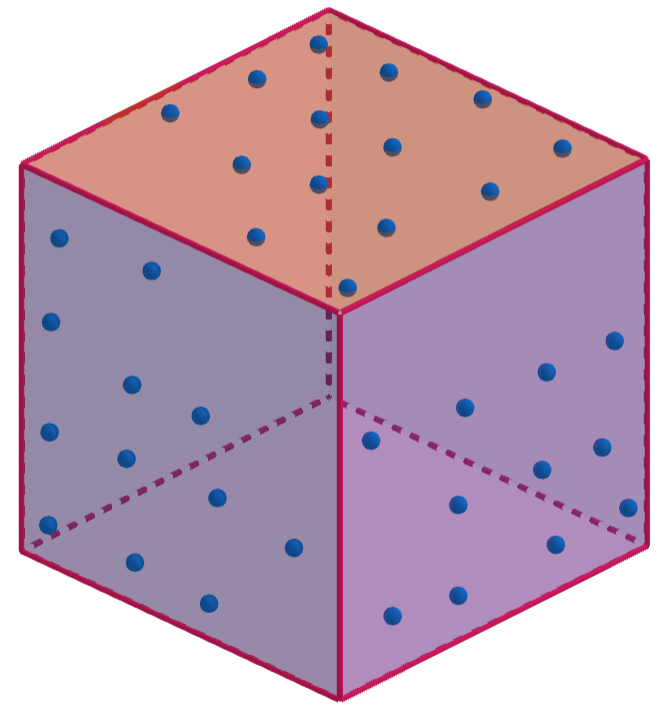}
\end{minipage}
}
\caption{\label{ss00} Simple random sampling.}
\end{figure*}

\begin{figure*}[h]
\centering
\subfigure[jittered sampling in two dimension]{
\begin{minipage}{7cm}
\centering
\includegraphics[width=0.6\textwidth]{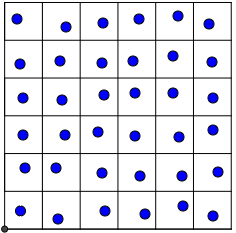}
\end{minipage}
}
\subfigure[jittered sampling in three dimension]{
\begin{minipage}{7cm}
\centering
\includegraphics[width=0.7\textwidth]{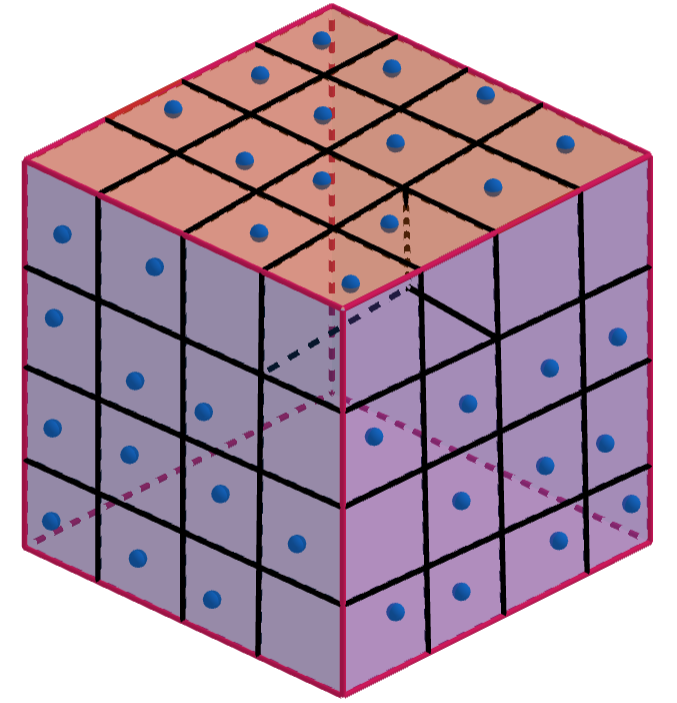}
\end{minipage}
}
\caption{\label{ss1} Jittered sampling formed by isometric grid partition.}
\end{figure*}

Stratified sampling is a special random sampling, that is different from \textbf{simple random sampling}, see Figure \ref{ss00}. The original sampling area is divided, and a uniformly distributed random sample point is selected in each subset of partitions. Jittered sampling is a special case of stratified sampling, involving grid-based equivolume partition. Explicitly, $[0,1]^{d}$ is divided into $m^d$ axis parallel boxes $Q_{i},1\leq i\leq N,$ each with sides $\frac{1}{m},$ see Figure \ref{ss1}. Research on the jittered sampling are extensive, see \cite{shirely1994,Doerr2,KP,KP2,jittsamp}.

We now consider a rectangle $R=[0,x)$ (we shall call it the test set in the following) in $[0,1]^{d}$ anchored at $0$. For an isometric grid partition $\Omega=\{Q_{1},Q_{2},
\ldots,Q_{N}\}$ of $[0,1]^{d}$, we put

$$
I_{N}:=\{j:\partial R \cap Q_{j}\neq \emptyset\},
$$
and

$$
C_{N}:=|I_{N}|,
$$
which means the cardinality of the index set $I_{N}$. For $C_N$, it is easy to obtain

\begin{equation}\label{CNbd0}
C_{N}\leq d\cdot N^{1-\frac{1}{d}}.
\end{equation}

\subsection{New partition models}\label{newpartmod}
In the end of this section, we design a class of partitions and construct it step by step. First, we consider the two-dimensional case.

\textbf{Step one: a class of partitions design for two dimension.}

\begin{figure}[H]
\centering
\includegraphics[width=0.30\textwidth]{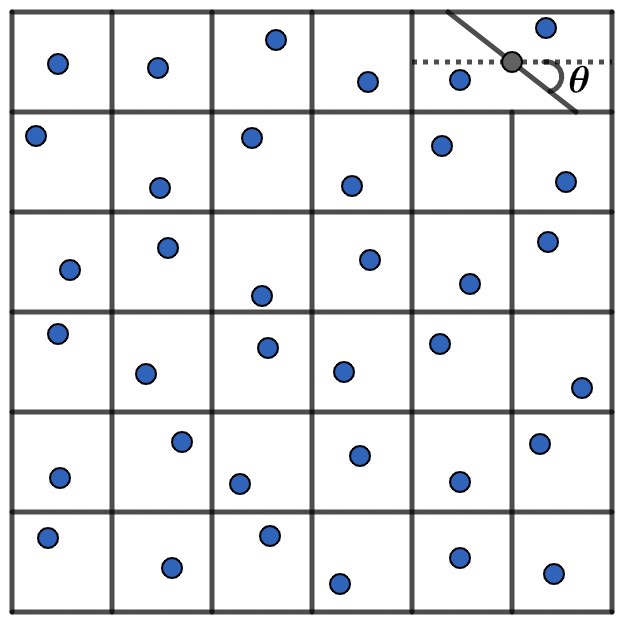}
\caption{\label{ss0}A class of partitions for two dimension}
\end{figure}

Our designed equivolume partition is actually a special case of general equivolume partition (see Figure \ref{ss0} for illustration in two dimensional case). For a grid-based equivolume partition in two dimension, we merge the two squares in the upper right corner to form a rectangle, then we use a series of straight line partitions to divide the rectangle into two equal-volume parts, which will be converted to a one-parameter model if we set the angle between the dividing line and horizontal line across the center $\theta$, where we suppose $0\leq\theta\leq\frac{\pi}{2}$. From simple calculations, we can conclude the arbitrary straight line must pass through \textbf{the center of the rectangle}. For convenience of notation, we set this partition model $\Omega_{\sim}=(\Omega_{1,\sim},\Omega_{2,\sim},Q_3,\ldots,Q_{N})$ in two dimensional case.

\begin{figure}[H]
\centering
\includegraphics[width=0.30\textwidth]{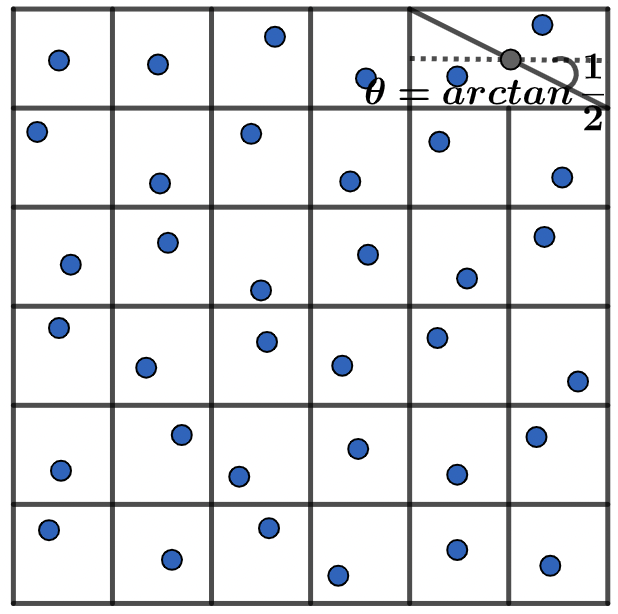}
\caption{\label{ss1111}The partition for parameter $\theta=arctan\ \frac{1}{2}$ in two dimension}
\end{figure}

In the above one-parameter model, the case will be grid-based equivolume partition if we choose $\theta=\frac{\pi}{2}$. The case $\theta=arctan\ \frac{1}{2}$ is introduced in \cite{KP}, see Figure \ref{ss1111} for two dimensional case. For notation convenience, we set this partition model $\Omega_{\backslash}=(\Omega_{1,\backslash},\Omega_{2,\backslash},Q_3,\ldots,Q_{N})$ in two dimensional case.

\begin{figure*}[h]
\centering
\subfigure[Isometric grid partition model]{
\begin{minipage}{7cm}
\centering
\includegraphics[width=1.0\textwidth]{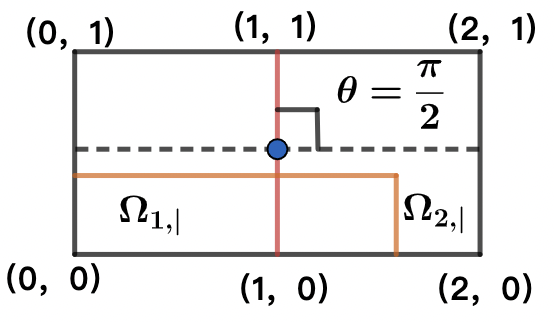}
\end{minipage}
}
\subfigure[Newly designed partition model]{
\begin{minipage}{7cm}
\centering
\includegraphics[width=1.0\textwidth]{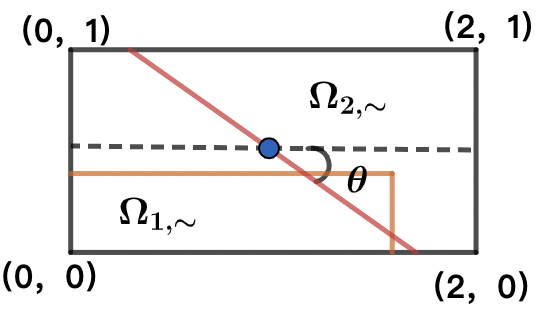}
\end{minipage}
}
\caption{\label{dbtpm} Difference between two partition models.}
\end{figure*}

The only difference between the new partition model and grid-based equivolume partition is to change two closed hypercubes into two special convex bodies, see illustration in Figure \ref{dbtpm}.

\textbf{Step two:} Suppose the original rectangle is $I$, for the convenience of calculation, we set the lower left corner of the rectangle at the origin $(0,0)$ and the side length of the small square to $1$. Now, consider $I=[0,2]\times[0,1]$ and its two equivolume partitions $(\Omega_{1,|},\Omega_{2,|})$ into two closed squares and $(\Omega_{1,\sim},\Omega_{2,\sim})$ into two convex bodies with

$$
    \Omega_{1,|}=[0,1]\times[0,1], \Omega_{1,\sim}=conv\{(0,0),(1+\frac{cot\theta}{2},0),(0,1),(1-\frac{cot\theta}{2},1)\},
$$
where $conv$ denotes the convex hull.


\textbf{Step three:} We consider the translation and stretch of the rectangle $I=[0,2]\times[0,1]$ into $$I'=[a_1,a_1+2b]\times [a_2,a_2+b],$$ the above two dimensional case in Step one can then be extended to $d-$dimension as \cite{KP}. Consider $d-$dimensional cuboid
\begin{equation}\label{eq28}
    I_d=I'\times\prod_{i=3}^{d}[a_i,a_i+b]
\end{equation} and its three equivolume partitions $\Omega'_{|}=(\Omega'_{1,|},\Omega'_{2,|})$ into two closed hypercubes, $\Omega'_{\backslash}=(\Omega'_{1,\backslash},\Omega'_{2,\backslash})$ into two closed, regular triangular hyperprisms and $\Omega'_{\sim}=(\Omega'_{1,\sim},\\ \Omega'_{2,\sim})$ into two closed, trapezoidal superconvex bodies with

\begin{equation}
    \Omega'_{1,|}=\prod_{i=1}^{d}[a_i,a_i+b],
\end{equation}

\begin{equation}\label{o1blaksla}
    \Omega'_{1,\backslash}=conv\{(a_1,a_2),(a_1+2b,a_2),(a_1,a_2+b)\}\times \prod_{i=3}^{d}[a_i,a_i+b],
\end{equation}
and

$$
\begin{aligned}
    \Omega'_{1,\sim}&=conv\{(a_1,a_2),(a_1+b+\frac{b\cdot cot\theta}{2},a_2),(a_1,a_2+b),(a_1+b-\frac{b\cdot cot\theta}{2},a_2+b)\}\\&\times \prod_{i=3}^{d}[a_i,a_i+b],
\end{aligned}
$$
where $conv$ denotes the convex hull. 

Just as grid-based partition $\bold{N=m^{d}}$, where $m$ represents the number of partitions in each dimension and $d$ denotes the dimensions. If we choose $a_1=\frac{m-2}{m}, a_2=\frac{m-1}{m}, b=\frac{1}{m}$, then, through the construction method from step one to step three, we get a
series of partitions (where we set $0\leq \theta\leq \frac{\pi}{2}$) that we call \textbf{local convex partition}, denoted by

\begin{equation}\label{omgasim}
\Omega^{*}_{\sim}=(\Omega^{*}_{1,\sim},\Omega^{*}_{2,\sim},Q_3 \ldots,Q_{N}).
\end{equation}

\begin{figure}[H]
\centering
\includegraphics[width=0.30\textwidth]{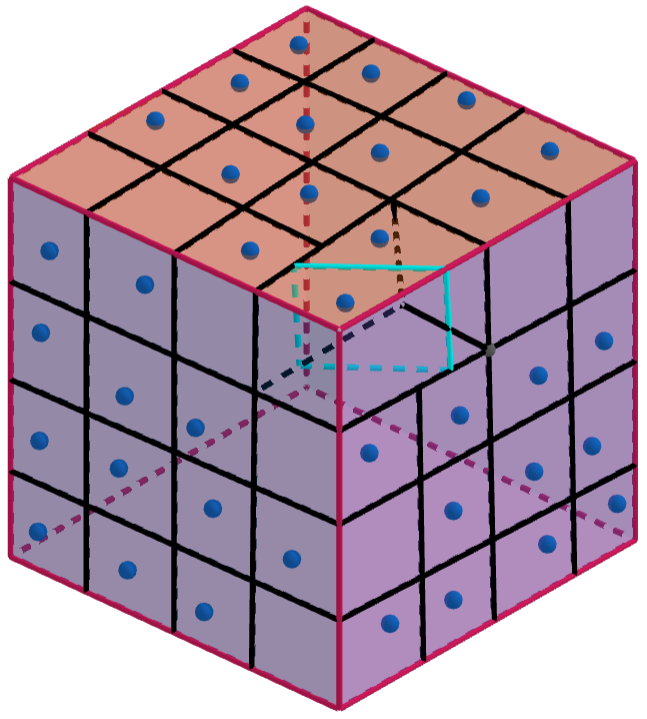}
\caption{Local convex partition in three dimension}
\end{figure}

Among the above local convex partition $\Omega^{*}_{\sim}$, if we choose the partition parameter $\theta=\frac{\pi}{2}$, isometric grid with partition number $m$ in each dimension is obtained, which we set
\begin{equation}\label{omga1}
\Omega^{*}_{|}=(\Omega^{*}_{1,|},\Omega^{*}_{2,|},Q_3 \ldots,Q_{N}).
\end{equation}

\begin{figure}[H]
\centering
\includegraphics[width=0.30\textwidth]{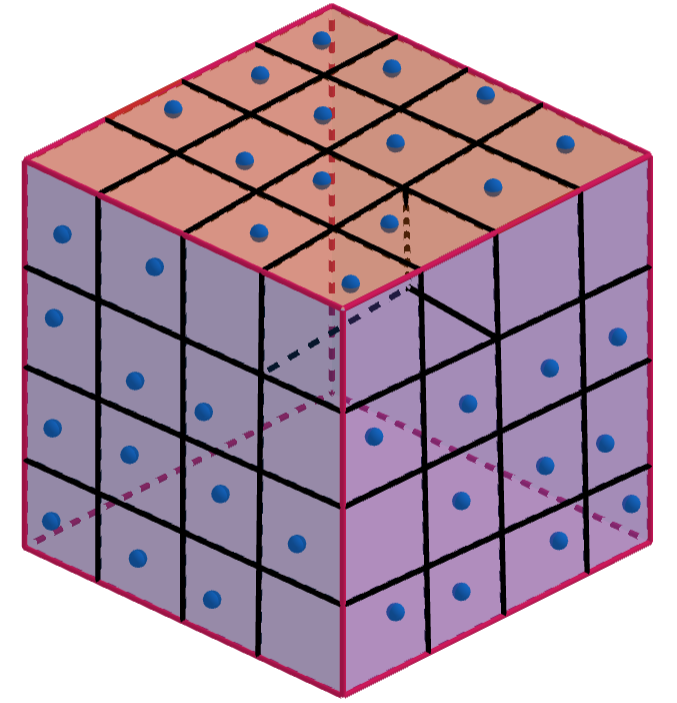}
\caption{local convex partition for parameter $\theta=\frac{\pi}{2}$ in three dimension}
\end{figure}

Likewise, if we choose the partition parameter $\theta=arctan\ \frac{1}{2}$, partition model in two dimensional case introduced above can then be extended to $d$ dimension, and we choose $a_1=\frac{m-2}{m}, a_2=\frac{m-1}{m}, b=\frac{1}{m}$ in \eqref{o1blaksla}, then this partition model is denoted by
\begin{equation}\label{omga2}
\Omega^{*}_{\backslash}=(\Omega^{*}_{1,\backslash},\Omega^{*}_{2,\backslash},Q_3 \ldots,Q_{N}).
\end{equation}

\begin{figure}[H]
\centering
\includegraphics[width=0.30\textwidth]{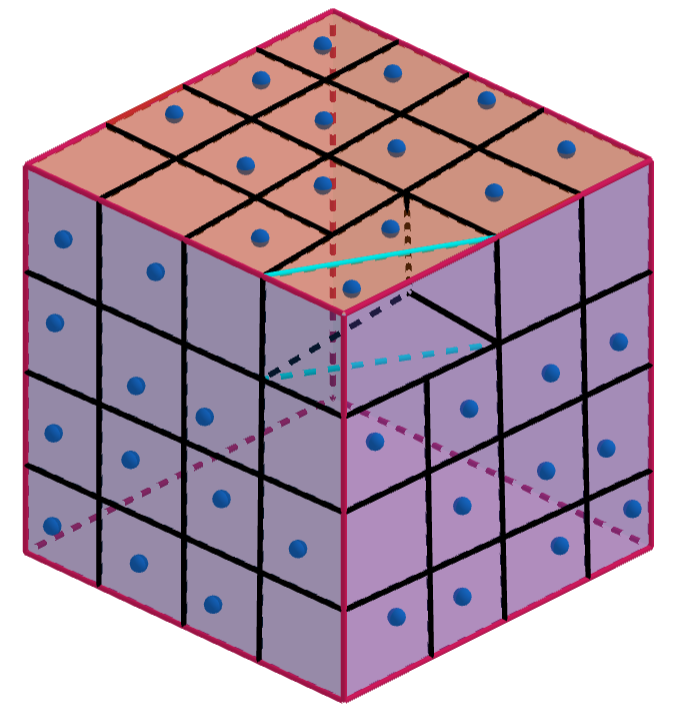}
\caption{local convex partition for parameter $\theta=arctan\ \frac{1}{2}$ in three dimension}
\end{figure}

\section{Expected $L_2-$discrepancy for stratified random sampling}\label{sec3}

In this section, comparisons of expected $L_2-$discrepancy under different partition models are obtained. Furthermore, we study expected $L_2-$discrepancy and several bounds are given under newly designed partition models.

\subsection{Expected $L_2-$discrepancy under two partition models}

\begin{thm}\label{uniformnoise0}
Let $m,d\in \mathbb{N}$ with $m\ge d\ge 2, 0\leq\theta\leq\frac{\pi}{2}$ and $N=m^d$. Stratified random $d-$dimension point sets $P_{\Omega^{*}_{|}}$ and $P_{\Omega^{*}_{\sim}}$ are uniformly distributed in the grid-based stratified subsets of $\Omega^{*}_{|}$ and stratified subsets of $\Omega^{*}_{\sim}$ respectively, then 

\begin{equation}\label{eq31}
\mathbb{E}(L_2^2(D_N,P_{\Omega^{*}_{\sim}}))\leq\mathbb{E}(L_2^2(D_N,P_{\Omega^{*}_{|}})).
\end{equation}
where $\Omega^{*}_{\sim}$, $\Omega^{*}_{|}$ are defined in \eqref{omgasim}, \eqref{omga1} respectively and $\theta$ is the partition parameter related to $\Omega^{*}_{\sim}$ as defined in Section $2$.
\end{thm}

\begin{rem}
In Theorem \ref{uniformnoise0}, as an infinite family of partitions is designed to generate point sets with a smaller expected $L_2-$discrepancy than classical stratified sampling (jittered sampling) for the same sampling number $N=m^d$. The equal signs on both sides of $\eqref{eq31}$ hold if and only if when $\theta=0$ or $\theta=\frac{\pi}{2}$.
\end{rem}

\begin{cor}\label{cor37}
Let  $N=m^d$ and $m,d\in \mathbb{N}$ with $m\ge d\ge 2$.  Stratified random $d-$dimension point sets $P_{\Omega^{*}_{|}}$ and $P_{\Omega^{*}_{\backslash}}$ are uniformly distributed in  $\Omega^{*}_{|}$ and $\Omega^{*}_{\backslash}$ respectively, then 

\begin{equation}\label{eq31}
\mathbb{E}(L_2^2(D_N,P_{\Omega^{*}_{\backslash}})) < \mathbb{E}(L_2^2(D_N,P_{\Omega^{*}_{|}})).
\end{equation}
where $\Omega^{*}_{|}$ and $\Omega^{*}_{\backslash}$ defined in  \eqref{omga1} and \eqref{omga2} respectively.
\end{cor}

\begin{rem}
Actually, (\ref{eq31}) holds if we choose parameter $\theta= arctan \frac{1}{2}$ in Theorem \ref{uniformnoise0}. The Corollary \ref{cor37} is main result in \cite{KP}. Obvious, the partition manner in \cite{KP} as Figure $4$ is included in our new partition models as Figure $3$.
\end{rem}

\subsection{Expected $L_2-$ discrepancy upper bounds under the new partition models}

In this subsection, expected $L_2-$discrepancy bounds under new partition models are given. Optimal result is also obtained under this class of partitions. 

\begin{thm}\label{expstarnp}
Let $m,d\in \mathbb{N}$ with $m\ge d\ge 2, 0\leq\theta\leq\frac{\pi}{2}$. Let $N=m^d$, the stratified random $d-$dimension point set $P_{\Omega^{*}_{\sim}}$ distributed in subsets of $\Omega^{*}_{\sim}$ defined in \eqref{omgasim}, then

\begin{equation}\label{ednstarx}
\mathbb{E}(L_2^2(D_N,P_{\Omega^{*}_{\sim}}))\leq \frac{d}{N^{1+\frac{1}{d}}}+\frac{1}{N^3}\cdot\frac{1}{3^{d-2}}\cdot P(\theta),
\end{equation} 
where

\begin{equation}\label{ptheta}
     P(\theta)=\left\{
\begin{aligned}
&\frac{2}{45}tan^3\theta+\frac{2}{15}tan^2\theta-\frac{tan\theta}{6}, \quad 0\leq\theta< arctan\frac{1}{2},
\\&-\frac{2}{45}, \quad  \theta=arctan\frac{1}{2},\\& -\frac{1}{24tan\theta}+\frac{1}{120tan^2\theta}+\frac{1}{1440tan^3\theta}, \quad  arctan\frac{1}{2}<\theta\leq\frac{\pi}{2}.
\end{aligned}
\right.
\end{equation}
\end{thm}

\begin{figure}[H]
\centering
\includegraphics[width=0.50\textwidth]{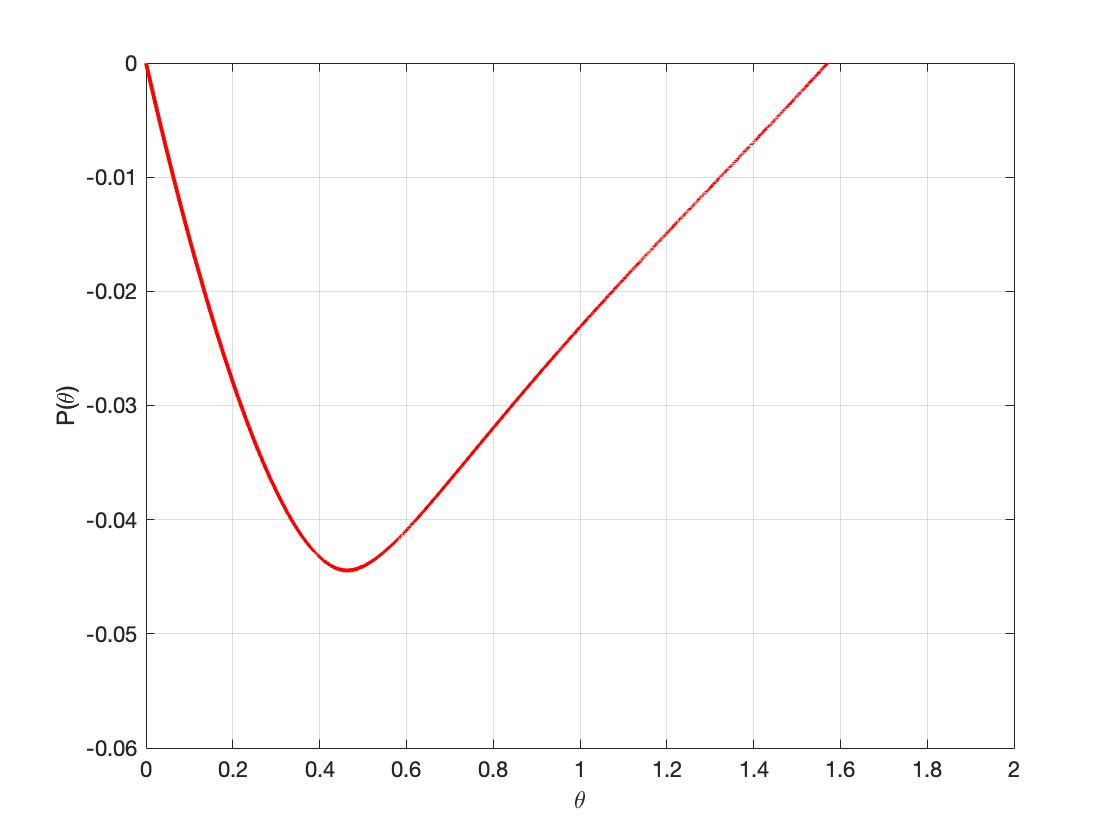}
\caption{\label{6p}$P(\theta)$ function}
\end{figure}

\begin{rem}
Noticing that in Theorem \ref{expstarnp}, $P(\theta)$ is a continuous function, decreases monotonically between $0$ and $arctan\frac{1}{2}$ and increases monotonically between $arctan\frac{1}{2}$ and $\frac{\pi}{2}$, see Figure \ref{6p}.  Choose parameter $\theta=\frac{\pi}{2}$ in Theorem \ref{expstarnp}, then we are back to the case of classical jittered sampling. Furthermore, all of these local convex partitions with parameter $\theta\in(0,\frac{\pi}{2})$ obtain better upper bounds of expected $L_2-$discrepancy than the jittered sampling. 
\end{rem}

\begin{cor}
Let $m,d\in \mathbb{N}$ with $m\ge d\ge 2$. Let $N=m^d$, the stratified random $d-$dimension point set $P_{\Omega^{*}_{\backslash}}$ distributed in subsets of $\Omega^{*}_{\backslash}$ defined in \eqref{omga2}, then we obtain optimal expected $L_2-$discrepancy bound under new partition models

\begin{equation}\label{ednstarx}
\mathbb{E}(L_2^2(D_N,P_{\Omega^{*}_{\backslash}}))\leq \frac{d}{N^{1+\frac{1}{d}}}-\frac{2}{45}\cdot\frac{1}{N^3}\cdot\frac{1}{3^{d-2}}.
\end{equation} 
\end{cor}

\begin{rem}
The optimal expected $L_2-$discrepancy bound under this class of partitions is obtained at $\theta=arctan\frac{1}{2}$ in Theorem \ref{expstarnp}. An upper bound on the expected $L_p-$discrepancy is derived by acceptance-rejection sampler using stratified inputs under the implicit constants in \cite{ZD2016}. Our results give \textbf{explicit expected $L_2-$discrepancy bounds} under a class of new  partitions,  which our order is the same with \cite{ZD2016}.
\end{rem}

\subsection{Some Examples}
This subsection presents some examples of expected $L_2-$\\discrepancy bounds under different sampling models for $N=m^d$. The cases of $\theta = arctan\ \frac{1}{2}$ and $\theta=\frac{\pi}{4}$ acquire better result than that of jittered sampling.

\textbf{Example 1.} Expected bound of stratified sampling set for $\theta=0$

$$\mathbb{E}(L_2^2(D_N))\leq \frac{d}{N^{1+\frac{1}{d}}}.$$

\begin{figure}[H]
\centering
\includegraphics[width=0.30\textwidth]{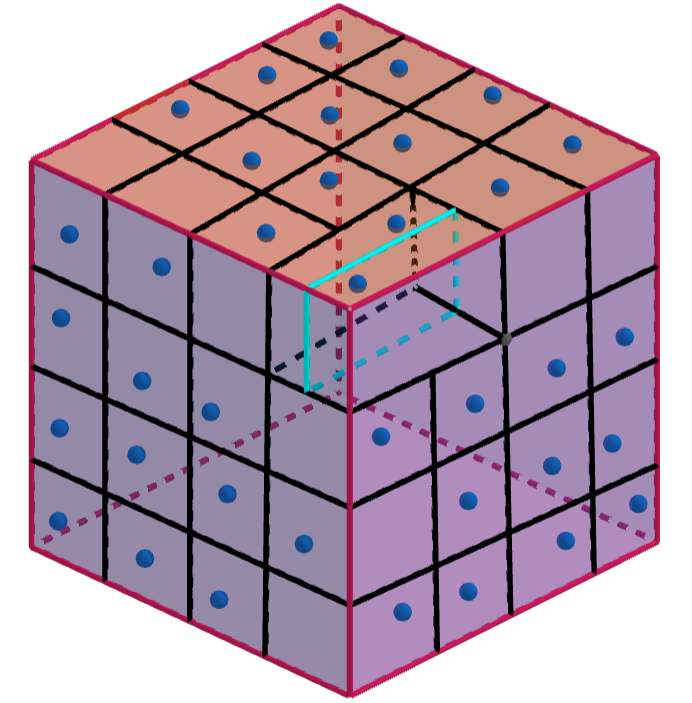}
\caption{Stratified sampling for $\theta=0$}
\end{figure}

\textbf{Example 2.} Expected bound of stratified sampling set for $\theta=\frac{\pi}{4}$

$$\mathbb{E}(L_2^2(D_N))\leq \frac{d}{N^{1+\frac{1}{d}}}-\frac{47}{1440\cdot 3^{d-2}\cdot N^3}.$$

\begin{figure}[H]
\centering
\includegraphics[width=0.30\textwidth]{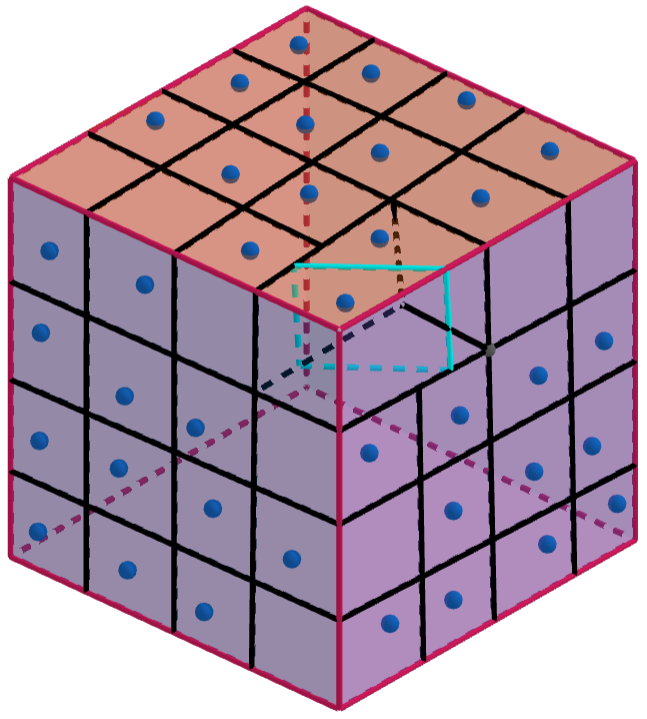}
\caption{Stratified sampling for $\theta=\frac{\pi}{4}$}
\end{figure}

\textbf{Example 3.} Expected bound of stratified sampling set for $\theta=\frac{\pi}{2}$

$$\mathbb{E}(L_2^2(D_N))\leq \frac{d}{N^{1+\frac{1}{d}}}.$$

\begin{figure}[H]
\centering
\includegraphics[width=0.30\textwidth]{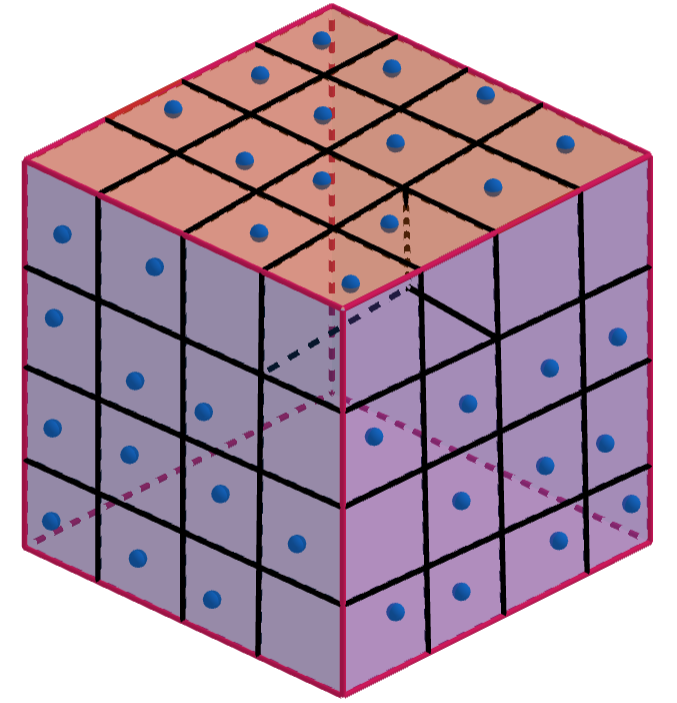}
\caption{Stratified sampling for $\theta=\frac{\pi}{2}$}
\end{figure}

\textbf{Example 4.} Expected bound of stratified sampling set for $\theta = arctan\ \frac{1}{2}$

$$\mathbb{E}(L_2^2(D_N))\leq \frac{d}{N^{1+\frac{1}{d}}}-\frac{2}{45\cdot 3^{d-2}\cdot N^3}.$$

\begin{figure}[H]
\centering
\includegraphics[width=0.30\textwidth]{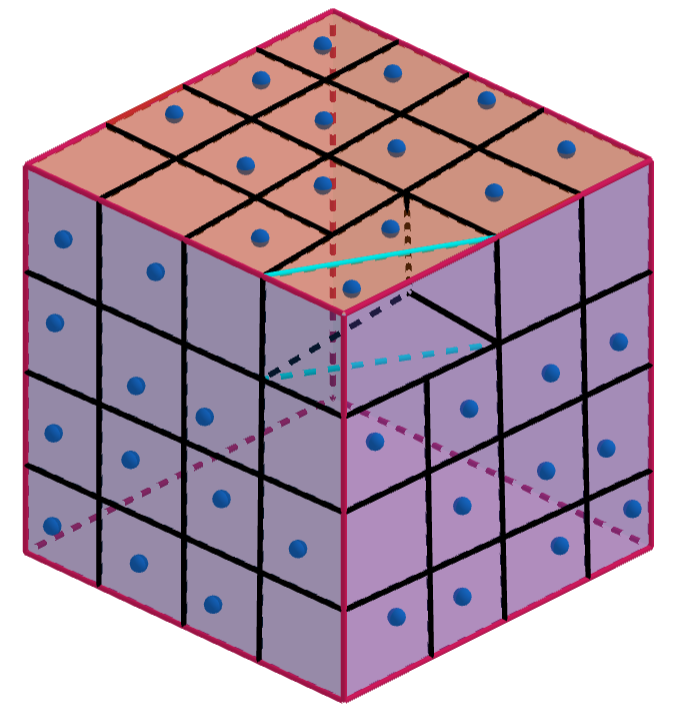}
\caption{Stratified sampling for $\theta = arctan\ \frac{1}{2}$}
\end{figure}

\section{Proofs}\label{pfmr}
In this section, we present the proofs of Theorem  \ref{uniformnoise0} and \ref{expstarnp}. The following lemma reveals the expected $L_2$-discrepancy quantitative relationship between the two partition models $\Omega^{*}_{|}$ and $\Omega^{*}_{\sim}$.

\begin{lma}\label{lm43}
For two equivolume partitions $\Omega^{*}_{\sim}=(\Omega^{*}_{1,\sim},\Omega^{*}_{2,\sim},Q_3 \ldots,Q_{N})$ and $\Omega^{*}_{|}=\{Q_1,Q_2,Q_3,\ldots,Q_N\}$ as defined in \eqref{omgasim} and \eqref{omga1} respectively, we have

\begin{equation}\label{ptheta}
     \mathbb{E}L_2^2(D_N,P_{\Omega^{*}_{\sim}})-\mathbb{E}L_2^2(D_N,P_{\Omega^{*}_{|}})=\left\{
\begin{aligned}
&\frac{1}{N^3}\cdot\frac{1}{3^{d-2}}\cdot P_{1}(\theta), \quad 0\leq\theta< arctan\frac{1}{2},
\\& -\frac{2}{45}\cdot \frac{1}{N^3}\cdot\frac{1}{3^{d-2}}, \quad \theta=arctan\frac{1}{2},\\& \frac{1}{N^3}\cdot\frac{1}{3^{d-2}}\cdot P_{2}(\theta), \quad  arctan\frac{1}{2}<\theta\leq\frac{\pi}{2}.
\end{aligned}
\right.
\end{equation}
where 

$$P_{1}(\theta)=\frac{2}{45}tan^3\theta+\frac{2}{15}tan^2\theta-\frac{tan\theta}{6},$$
and

$$P_{2}(\theta)=-\frac{1}{24tan\theta}+\frac{1}{120tan^2\theta}+\frac{1}{1440tan^3\theta}.$$
\end{lma}

\subsection{Proof of Lemma \ref{lm43}}

For equivolume partition $\mathbf{\Omega_{0,\sim}}=(\Omega_{1,\sim},\Omega_{2,\sim})$ of $I$(the same argument if we replace $\mathbf{\Omega_{0,\sim}}$ with $\mathbf{\Omega_{0,|}}$), from [Proposition $2$] in \cite{KP}, which is, for an equivolume partition $\Omega=\{\Omega_1, \Omega_2, \ldots, \Omega_N\}$ of a compact convex set $K\subset \mathbb{R}^{d}$ with $\lambda(K)>0$, $P_{\Omega}$ is the corresponding stratified sampling set, then 

\begin{equation}\label{el22pom}
    \mathbb{E}L_2^2(D_N, P_{\Omega})=\frac{1}{N^2\lambda(K)}\sum_{i=1}^{N}\int_{K}q_i(x)(1-q_i(x))dx,
\end{equation}
where 

\begin{equation}\label{qixx}
    q_i(x)=\frac{\lambda(\Omega_i\cap[0,x])}{\lambda(\Omega_i)}.
\end{equation}

Through simple derivation, it follows that
 \begin{equation}\label{18qi}
    \mathbb{E}L_2^2(D_N,P_{\mathbf{\Omega_{0,\sim}}})=\frac{1}{8}\sum_{i=1}^{2}\int_{I}\mathbf{q}_i(x)(1-\mathbf{q}_i(x))dx,
 \end{equation}
and

\begin{equation}\label{qix}
    \mathbf{q}_i(x)=\frac{\lambda(\Omega_{i,\sim}\cap [0,x])}{\lambda(\Omega_{i,\sim})}=\lambda(\Omega_{i,\sim}\cap [0,x]).
\end{equation}

Conclusion \eqref{18qi} is equivalent to the following

$$
    8\mathbb{E}L_2^2(D_N,P_{\mathbf{\Omega_{0,\sim}}})=1-\sum_{i=1}^{2}\int_{I}\mathbf{q}_i^{2}(x)dx.
$$

We first consider parameter $arctan\frac{1}{2}\leq\theta\leq\frac{\pi}{2}$, then we define the following two functions for simplicity of the expression.

$$
    F(\mathbf{x})=\frac{1}{2}\cdot[(x_1-1)tan\theta+x_2-\frac{1}{2}]\cdot[(x_1-1)+(x_2-\frac{1}{2})\cdot cot\theta],
$$
and

$$
    G(\mathbf{x})=x_1x_2-x_2-\frac{cot\theta}{2}x_2+\frac{1}{2}x_2^2\cdot cot\theta,
$$
where $\mathbf{x}=(x_1,x_2)$.

Furthermore, for $\mathbf{\Omega_{0,|}}=(\Omega_{1,|},\Omega_{2,|})$ defined in Step two of Section \ref{newpartmod}, \eqref{qix} implies

$$
    \mathbf{q}_{1,|}(\mathbf{x})=\left\{
\begin{aligned}
&x_1x_2, \mathbf{x}\in \Omega_{1,|}\\&
x_2, \mathbf{x}\in \Omega_{2,|},
\end{aligned}
\right.
$$
and 
$$
    \mathbf{q}_{2,|}(\mathbf{x})=\left\{
\begin{aligned}
&0, \mathbf{x}\in \Omega_{1,|}\\&
(x_1-1)x_2, \mathbf{x}\in \Omega_{2,|}.
\end{aligned}
\right.
$$

\begin{figure}[H]
\centering
\includegraphics[width=0.60\textwidth]{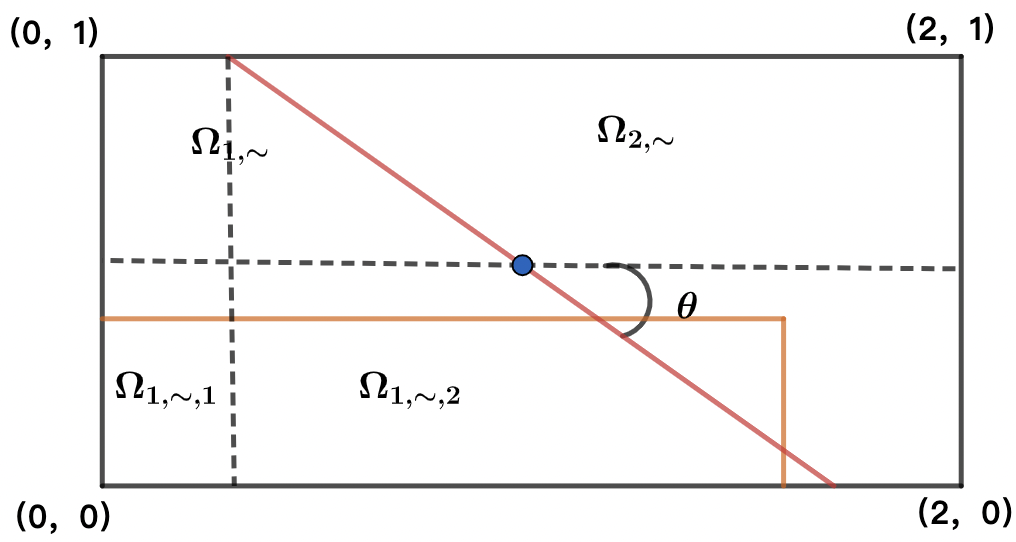}
\caption{Division of the integral region}\label{spt1}
\end{figure}

\begin{figure*}[h]
\centering
\subfigure[]{
\begin{minipage}{7cm}
\centering
\includegraphics[width=1.0\textwidth]{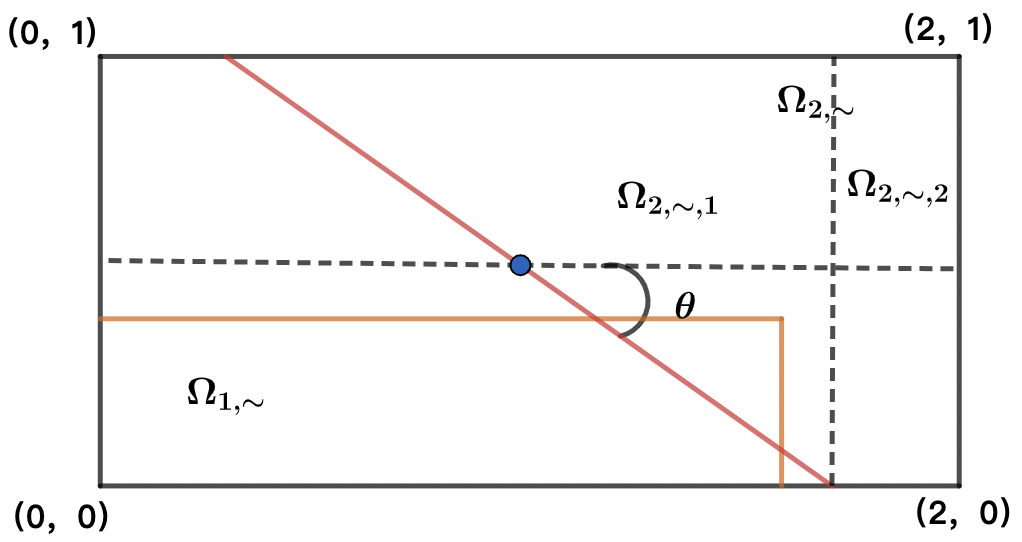}
\end{minipage}
}
\subfigure[]{
\begin{minipage}{7cm}
\centering
\includegraphics[width=1.0\textwidth]{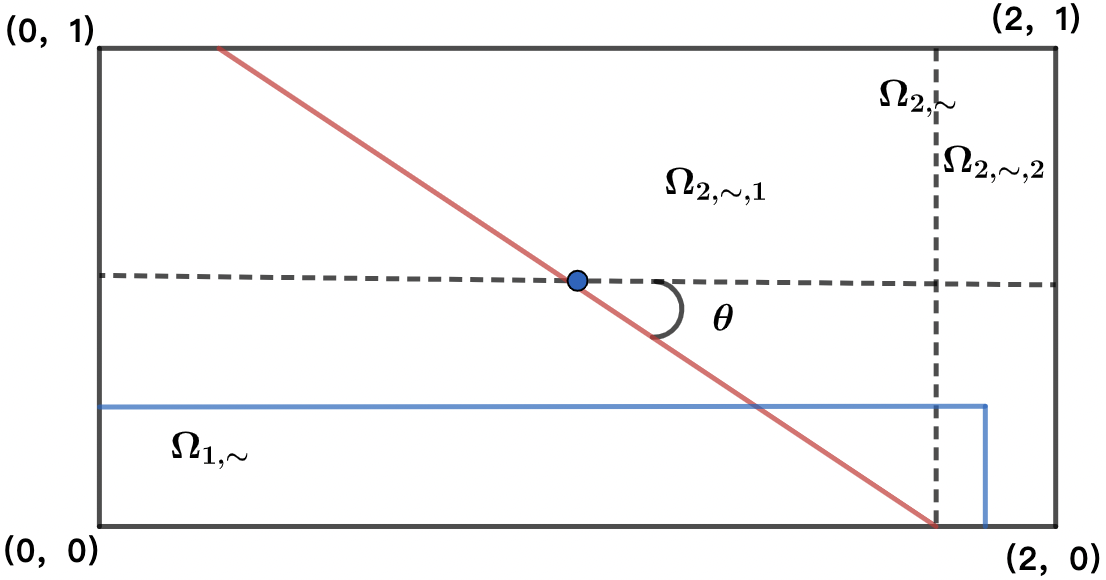}
\end{minipage}
}
\caption{\label{figure6} Division of the integral region.}
\end{figure*}

Besides, 

$$
    \mathbf{q}_{1,\sim}(\mathbf{x})=\left\{
\begin{aligned}
&x_1x_2, \mathbf{x}\in \Omega_{1,\sim},\\&
x_1x_2-F(\mathbf{x}), \mathbf{x}\in \Omega_{2,\sim,1},\\&x_1x_2-G(\mathbf{x}), \mathbf{x}\in \Omega_{2,\sim,2},
\end{aligned}
\right.
$$
and 
$$
    \mathbf{q}_{2,\sim}(\mathbf{x})=\left\{
\begin{aligned}
&0, \mathbf{x}\in \Omega_{1,\sim},\\&
F(\mathbf{x}), \mathbf{x}\in \Omega_{2,\sim,1},\\&G(\mathbf{x}), \mathbf{x}\in \Omega_{2,\sim,2},
\end{aligned}
\right.
$$
where $\Omega_{1,\sim}$, $\Omega_{2,\sim}$ denote subsets of partition $\mathbf{\Omega_{0,\sim}}$. In the following, we shall continue to divide subsets $\Omega_{1,\sim}=\{\Omega_{1,\sim,1},\Omega_{1,\sim,2}\}$ and $\Omega_{2,\sim}=\{\Omega_{2,\sim,1},\Omega_{2,\sim,2}\}$ to facilitate calculation. See Figures \ref{spt1} to \ref{figure6}.

Therefore, for $\theta=\frac{\pi}{2}$, we introduce two symbols $B_{1,|},B_{2,|}$ and have

\begin{equation}\label{b11}
    B_{1,|}=\int_{I}\mathbf{q}_{1,|}^2(\mathbf{x})d\mathbf{x}=\int_{\Omega_{1,|}}x_1^2x_2^2d\mathbf{x}+\int_{\Omega_{2,|}}x_2^2d\mathbf{x}=\frac{1}{9}+\frac{1}{3}=\frac{4}{9},
\end{equation}
and

\begin{equation}\label{b21}
  B_{2,|}=\int_{I}\mathbf{q}_{2,|}^2(\mathbf{x})d\mathbf{x}=\int_{\Omega_{2,|}}(x_1-1)^2x_2^2d\mathbf{x}=\frac{1}{9}.
\end{equation}

Thus,

\begin{equation}\label{8el2p1}
    8\mathbb{E}(L_2^2(P_{\mathbf{\Omega_{0,|}}}))=1-(B_{1,|}+B_{2,|})=\frac{4}{9}.
\end{equation}

Furthermore, we introduce $B_{1,\sim}$ and $B_{2,\sim}$, then

\begin{equation}\label{b1sim}
\begin{aligned}
    B_{1,\sim}=\int_{I}\mathbf{q}_{1,\sim}^2(\mathbf{x})d\mathbf{x}&=\int_{\Omega_{1,\sim}}x_1^2x_2^2d\mathbf{x}+\int_{\Omega_{2,\sim,1}}(x_1x_2-F(\mathbf{x}))^2d\mathbf{x}\\&+\int_{\Omega_{2,\sim,2}}(x_1x_2-G(\mathbf{x}))^2d\mathbf{x},
    \end{aligned}
\end{equation}
and

$$
  B_{2,\sim}=\int_{I}\mathbf{q}_{2,\sim}^2(\mathbf{x})d\mathbf{x}=\int_{\Omega_{2,\sim,1}}F^{2}(\mathbf{x})d\mathbf{x}+\int_{\Omega_{2,\sim,2}}G^{2}(\mathbf{x})d\mathbf{x}.$$
  
We divide our calculation in three steps. First, we compute $\int_{\Omega_{1,\sim}}x_1^2x_2^2d\mathbf{x}$, see Figure \ref{spt1} for illustration.

\begin{equation}\label{omega1sim1}
\begin{aligned}
    \int_{\Omega_{1,\sim,1}}x_1^2x_2^2d\mathbf{x}&=\int_{0}^{1-\frac{cot\theta}{2}}x_1^{2}dx_1\cdot \int_{0}^{1}x_2^{2}dx_2=\frac{(2-cot\theta)^{3}}{72}.
    \end{aligned}
\end{equation}

\begin{equation}\label{omega1sim2}
\begin{aligned}
    \int_{\Omega_{1,\sim,2}}x_1^2x_2^2d\mathbf{x}&=\int_{1-\frac{cot\theta}{2}}^{1+\frac{cot\theta}{2}}x_1^{2}dx_1\cdot \int_{0}^{(1-x_1)\cdot tan\theta+\frac{1}{2}}x_2^{2}dx_2
    \\&=\frac{60tan^{2}\theta-36tan\theta+7}{720tan^{3}\theta}.
    \end{aligned}
\end{equation}

Therefore, \eqref{omega1sim1} and \eqref{omega1sim2} imply

\begin{equation}\label{ome1sim}
\begin{aligned}
    \int_{\Omega_{1,\sim}}x_1^2x_2^2d\mathbf{x}&=\int_{\Omega_{1,\sim,1}}x_1^2x_2^2d\mathbf{x}+\int_{\Omega_{1,\sim,2}}x_1^2x_2^2d\mathbf{x}
    \\&=-\frac{1}{12tan\theta}+\frac{1}{30tan^2\theta}-\frac{1}{240tan^3\theta}+\frac{1}{9}.
\end{aligned}
\end{equation}

Second, we compute $\int_{\Omega_{2,\sim,1}}(x_1x_2-F(\mathbf{x}))^2d\mathbf{x}$ and $\int_{\Omega_{2,\sim,2}}(x_1x_2-G(\mathbf{x}))^2d\mathbf{x}$.

\begin{equation}\label{omega2sim1}
\begin{aligned}
    \int_{\Omega_{2,\sim,1}}(x_1x_2-F(\mathbf{x}))^2d\mathbf{x}&=\int_{1-\frac{cot\theta}{2}}^{1+\frac{cot\theta}{2}}\int_{(1-x_1)\cdot tan\theta+\frac{1}{2}}^{1}(x_1x_2-F(\mathbf{x}))^2dx_2dx_1\\&=\frac{180tan^{2}\theta-12tan\theta+5}{720tan^{3}\theta},
    \end{aligned}
\end{equation}

\begin{equation}\label{omega2sim2}
\begin{aligned}
    &\int_{\Omega_{2,\sim,2}}(x_1x_2-G(\mathbf{x}))^2d\mathbf{x}=\int_{1+\frac{cot\theta}{2}}^{2}\int_{0}^{1}(x_1x_2-G(\mathbf{x}))^2dx_2dx_1\\&=-\frac{cot^3\theta}{240}-\frac{cot^2\theta}{30}-\frac{cot\theta}{12}+\frac{1}{3}.
    \end{aligned}
\end{equation}

Thus, \eqref{omega2sim1} and \eqref{omega2sim2} imply 

\begin{equation}\label{ome2sim23}
\begin{aligned}
    &\int_{\Omega_{2,\sim,1}}(x_1x_2-F(\mathbf{x}))^2d\mathbf{x}+\int_{\Omega_{2,\sim,2}}(x_1x_2-G(\mathbf{x}))^2d\mathbf{x}\\&=\frac{1}{3}+\frac{1}{6tan\theta}-\frac{1}{20tan^2\theta}+\frac{1}{360tan^3\theta}.
\end{aligned}
\end{equation}

Combining \eqref{b1sim}, \eqref{ome1sim} and \eqref{ome2sim23}, we have

\begin{equation}\label{B1sim1}
\begin{aligned}
     B_{1,\sim}&=\int_{\Omega_{1,\sim}}x_1^2x_2^2d\mathbf{x}+\int_{\Omega_{2,\sim,1}}(x_1x_2-F(\mathbf{x}))^2d\mathbf{x}+\int_{\Omega_{2,\sim,2}}(x_1x_2-G(\mathbf{x}))^2d\mathbf{x}\\&=\frac{1}{12tan\theta}-\frac{1}{60tan^2\theta}-\frac{1}{720tan^3\theta}+\frac{4}{9}.
    \end{aligned}
\end{equation}

Third, we will compute $\int_{\Omega_{2,\sim,1}}F^{2}(\mathbf{x})d\mathbf{x}$ and $\int_{\Omega_{2,\sim,2}}G^2(\mathbf{x})d\mathbf{x}$ in the following.

In fact,

\begin{equation}\label{omega22sim1}
\begin{aligned}
    \int_{\Omega_{2,\sim,1}}F^{2}(\mathbf{x})d\mathbf{x}&=\int_{1-\frac{cot\theta}{2}}^{1+\frac{cot\theta}{2}}\int_{(1-x_1)\cdot tan\theta+\frac{1}{2}}^{1}F^{2}(\mathbf{x})dx_2dx_1\\&=\frac{1}{120tan^{3}\theta},
    \end{aligned}
\end{equation}

\begin{equation}\label{omega22sim2}
\begin{aligned}
    &\int_{\Omega_{2,\sim,2}}G^{2}(\mathbf{x})d\mathbf{x}=\int_{1+\frac{cot\theta}{2}}^{2}\int_{0}^{1}G^{2}(\mathbf{x})dx_2dx_1\\&=\frac{1}{9}-\frac{1}{24tan\theta}+\frac{1}{120tan^2\theta}-\frac{11}{1440tan^3\theta}.
    \end{aligned}
\end{equation}

Combining \eqref{omega22sim1} and \eqref{omega22sim2}, we have

\begin{equation}\label{ome2sim}
\begin{aligned}
    B_{2,\sim}&=\int_{\Omega_{2,\sim,1}}F^2(\mathbf{x})d\mathbf{x}+\int_{\Omega_{2,\sim,2}}G^2(\mathbf{x})d\mathbf{x}\\&=\frac{1}{9}-\frac{1}{24tan\theta}+\frac{1}{120tan^2\theta}+\frac{1}{1440tan^3\theta}.
\end{aligned}
\end{equation}

Thus, 

\begin{equation}\label{b1simb2sim}
   B_{1,\sim}+B_{2,\sim}=\frac{1}{24tan\theta}-\frac{1}{120tan^2\theta}-\frac{1}{1440tan^3\theta}+\frac{5}{9}.
\end{equation}

Therefore,

\begin{equation}\label{el2p}
\begin{aligned}
    8\mathbb{E}(L_2^2(P_{\mathbf{\Omega_{0,\sim}}}))&=1-(B_{1,\sim}+B_{2,\sim})
    \\&=-\frac{cot\theta}{24}+\frac{cot^2\theta}{120}+\frac{cot^3\theta}{1440}+\frac{4}{9},
    \end{aligned}
\end{equation}
where $arctan\frac{1}{2}\leq\theta<\frac{\pi}{2}$. 

For $\theta=\frac{\pi}{2}$, by \eqref{8el2p1} we have

\begin{equation}\label{el22149}
8\mathbb{E}(L_2^2(P_{\mathbf{\Omega_{0,|}}}))=\frac{4}{9}.
\end{equation}

\begin{figure*}[h]
\centering
\subfigure[]{
\begin{minipage}{7cm}
\centering
\includegraphics[width=1.0\textwidth]{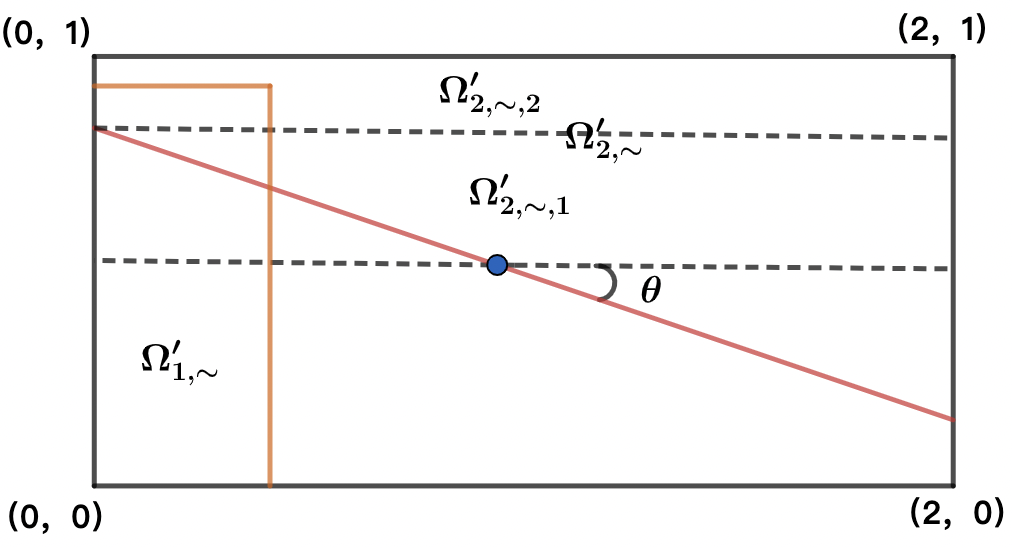}
\end{minipage}
}
\subfigure[]{
\begin{minipage}{7cm}
\centering
\includegraphics[width=1.0\textwidth]{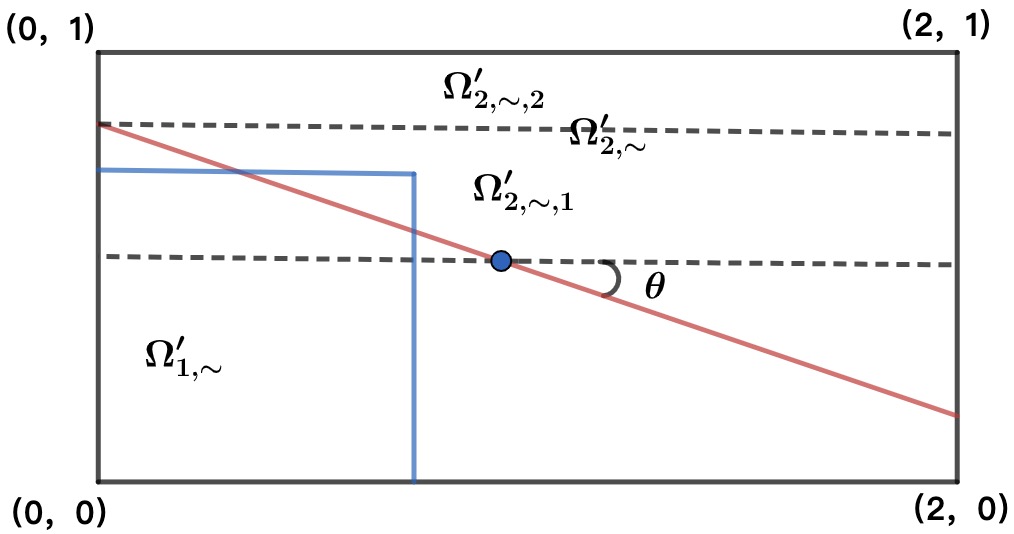}
\end{minipage}
}
\caption{\label{figure7} Division of the integral region.}
\end{figure*}

Considering the case $0\leq\theta<arctan\frac{1}{2}$, we denote the partition by $\Omega'_{\sim}=\{\Omega'_{1,\sim},\Omega'_{2,\sim}\}$, see Figure \ref{figure7}. Let

$$
    \mathbf{q}'_{1,\sim}(\mathbf{x})=\left\{
\begin{aligned}
&x_1x_2, \mathbf{x}\in \Omega'_{1,\sim},\\&
x_1x_2-H(\mathbf{x}), \mathbf{x}\in \Omega'_{2,\sim,1},\\& x_1x_2-J(\mathbf{x}), \mathbf{x}\in \Omega'_{2,\sim,2}.
\end{aligned}
\right.
$$
and 
$$
    \mathbf{q}'_{2,\sim}(\mathbf{x})=\left\{
\begin{aligned}
&0, \mathbf{x}\in \Omega'_{1,\sim},\\&
H(\mathbf{x}), \mathbf{x}\in \Omega'_{2,\sim,1},\\&
J(\mathbf{x}), \mathbf{x}\in \Omega'_{2,\sim,2},
\end{aligned}
\right.
$$
where
\begin{equation}
    H(x)=\frac{1}{2}\cdot[x_2-(1-x_1)tan\theta-\frac{1}{2}]\cdot[cot\theta\cdot x_2-1+x_1-\frac{1}{2}cot\theta],
\end{equation}
and

\begin{equation}
    J(x)=[x_2-tan\theta-\frac{1}{2}]\cdot x_1+\frac{1}{2}x_1^2\cdot tan\theta.
\end{equation}
Then we divide subsets $\Omega'_{1,\sim}=\{\Omega'_{1,\sim,1},\Omega'_{1,\sim,2}\}$ and $\Omega'_{2,\sim}=\{\Omega'_{2,\sim,1}, \Omega'_{2,\sim,2}\}$ to facilitate calculation. See Figure \ref{figure7}.

So 

\begin{equation}\label{eq424}
\begin{aligned}
    B'_{1,\sim}=\int_{I}\mathbf{q}_{1,\sim}^{'2}(\mathbf{x})d\mathbf{x}&=\int_{\Omega'_{1,\sim}}x_1^2x_2^2d\mathbf{x}+\int_{\Omega'_{2,\sim,1}}(x_1x_2-H(\mathbf{x}))^2d\mathbf{x}\\&+\int_{\Omega'_{2,\sim,2}}(x_1x_2-J(\mathbf{x}))^2d\mathbf{x},
    \end{aligned}
\end{equation}
and

$$
  B'_{2,\sim}=\int_{I}\mathbf{q}_{2,\sim}^{'2}(\mathbf{x})d\mathbf{x}=\int_{\Omega'_{2,\sim,1}}H^{2}(\mathbf{x})d\mathbf{x}+\int_{\Omega'_{2,\sim,2}}J^{2}(\mathbf{x})d\mathbf{x}.
$$

If we follow the calculation process of \eqref{omega1sim1}-\eqref{b1simb2sim}, then we obtain 

\begin{equation}
    B'_{1,\sim}=-\frac{4}{45}tan^3\theta-\frac{4}{15}tan^2\theta+\frac{tan\theta}{3}+\frac{4}{9},
\end{equation}
and

\begin{equation}\label{b2sim}
    B'_{2,\sim}=\frac{2}{45}tan^3\theta+\frac{2}{15}tan^2\theta-\frac{tan\theta}{6}+\frac{1}{9}.
\end{equation}

Thus,

\begin{equation}\label{b1simb2sim111}
   B'_{1,\sim}+B'_{2,\sim}=-\frac{2}{45}tan^3\theta-\frac{2}{15}tan^2\theta+\frac{tan\theta}{6}+\frac{5}{9}.
\end{equation}

Hence,

\begin{equation}\label{lygjd}
\begin{aligned}
    8\mathbb{E}(L_2^2(P_{\Omega'_{\sim}}))&=1-(B'_{1,\sim}+ B'_{2,\sim})\\&=\frac{4}{9}+\frac{2}{45}tan^3\theta+\frac{2}{15}tan^2\theta-\frac{tan\theta}{6},
    \end{aligned}
\end{equation}
where $0\leq\theta<arctan\frac{1}{2}$.

Combining with \eqref{el2p} and considering the translation and stretch of the rectangle $I=[0,2]\times[0,1]$ into $$I'=[a_1,a_1+2b]\times [a_2,a_2+b],$$ we obtain

\begin{equation}\label{440comp28}
\mathbb{E}(L_2^2(P_{\Omega^{*}_{\sim}}))\leq\mathbb{E}(L_2^2(P_{\Omega^{*}_{|}})),
\end{equation}
where $a_1=\frac{m-2}{m},a_2=\frac{m-1}{m},b=\frac{1}{m}$, $\Omega^{*}_{\sim}$ is the infinite family of equivolume partitions defined in \eqref{omgasim} and $\Omega^{*}_{|}$ is grid-based equivolume partition defined in \eqref{omga1}. The equal sign of \eqref{440comp28} holds if and only if partition parameter $\theta=0,\frac{\pi}{2}$. Noting that conclusion \eqref{440comp28} is only for \textbf{the two-dimensional case}.

\textbf{Next we will give a proof of \eqref{440comp28} for $d-$dimensional case.} We firstly prove the case $b=1$ and $(a_1,a_2,\ldots,a_d)=(0,0,\ldots,0).$ Let $I'_d=[0,2]\times[0,1]\times[0,1]^{d-2}$ and we denote partition manner of this special case $\Omega''_{\sim}=\{\Omega''_{1,\sim},\Omega''_{2,\sim}\}$.  

For $i=1,2$, we have

$$
    \mathbf{q}'_{i,\sim}(\mathbf{x})=\mathbf{q}_{i,\sim}(x_1,x_2)\cdot\prod_{j=3}^{d}x_j,
$$
where $\mathbf{q}'_{i,\sim}(\mathbf{x})$ is defined as \eqref{qix} for $\Omega''_{\sim}$.

Thus,

$$
    \int_{I'_d} \mathbf{q}'^{2}_{i,\sim}(\mathbf{x})d\mathbf{x}=B_{i,\sim}\cdot\int_{[0,1]^{d-2}}\prod_{j=3}^{d}x_j^2dx_3dx_4\ldots dx_d=\frac{1}{3^{d-2}}\cdot B_{i,\sim},
$$
where $B_{i,\sim},i=1,2$ have been calculated in \eqref{B1sim1} and \eqref{ome2sim} respectively.

As we have 

$$
    \int_{I'_d}\lambda([0,\mathbf{x}])d\mathbf{x}=\int_{[0,1]^{d-2}}\prod_{j=3}^{d}x_jdx_3dx_4\ldots dx_d=\frac{1}{2^{d-2}}.
$$

Then we obtain,

\begin{equation}\label{8l22omesim}
8\mathbb{E}(L_2^2(P_{\Omega''_{\sim}}))=\frac{1}{2^{d-2}}-\frac{1}{3^{d-2}}\cdot(B_{1,\sim}+B_{2,\sim}).
\end{equation}

Now, for $I_d$ in \eqref{eq28}, we define a vector

\begin{equation}
    \mathbf{a}=\{a_1,a_2,\ldots,a_d\}.
\end{equation}

We then prove \eqref{el22pom} is independent of $\mathbf{a}$. In $I_d$, we choose $\mathbf{a}=0$, set

\begin{equation}\label{id0433}
    I_{d}^{0}=[0,2b]\times [0,b]^{d-1},
\end{equation}
and
\begin{equation}\label{id0m433}
    I_{d,m}^{0}=[0,\frac{2}{m}]\times [0,\frac{1}{m}]^{d-1}.
\end{equation}

It suffices to show that

\begin{equation}\label{eq434}
     \frac{1}{N^2\lambda(I_d)}\sum_{i=1}^{N}\int_{I_d}q_i(x)(1-q_i(x))dx=\frac{1}{N^2\lambda(I_{d}^{0})}\sum_{i=1}^{N}\int_{I_{d}^{0}}q_i(x)(1-q_i(x))dx.
\end{equation}

We only consider $N=2$ in \eqref{eq434}, this is because we choose $K=I_d$ and $K=I_{d}^{0}$ in \eqref{el22pom} respectively. This means $I_d,I_{d}^{0}$ are divided into two equal volume parts respectively.

Let 

\begin{equation}\label{xiaiti}
    x_i-a_i=t_i,1\leq i\leq d.
\end{equation}

According to \eqref{qixx} and plugging \eqref{xiaiti} into the left side of \eqref{eq434}, the desired result is obtained.

From \eqref{el22pom} and let $K=[0,1]^d$, we have

\begin{equation}\label{eq437}
\begin{aligned}
    &\mathbb{E}L_2^2(P_{\Omega^{*}_{\sim}})-\mathbb{E}L_2^2(P_{\Omega^{*}_{|}})\\&=\frac{1}{N^2}\sum_{i=1}^{N}\int_{[0,1]^d}\tilde{q}_i(x)(1-\tilde{q}_i(x))dx-\frac{1}{N^2}\sum_{i=1}^{N}\int_{[0,1]^d}\bar{q}_i(x)(1-\bar{q}_i(x))dx,
\end{aligned}
\end{equation}
where 
$$ \tilde{q}_i(x)=\frac{\lambda(\Omega^{*}_{i,\sim}\cap[0,x])}{\lambda(\Omega^{*}_{i,\sim})}, \bar{q}_i(x)=\frac{\lambda(\Omega^{*}_{i,|}\cap[0,x])}{\lambda(\Omega^{*}_{i,|})}, i=1,2,$$
and 
$$\tilde{q}_i(x)=\bar{q}_i(x)=\frac{\lambda(Q_i\cap[0,x])}{\lambda(Q_i)}, i=3,4,\ldots,N.$$

Let $I_{d,m}^{0}=\{\Omega^{*}_{1,\sim},\Omega^{*}_{2,\sim}\}$, $I_{d,m}^{0}=\{\Omega^{*}_{1,|},\Omega^{*}_{2,|}\}$ denote two different partitions of $I_{d,m}^{0}$. It can easily be seen only $I_{d,m}^{0}$ contributes to the difference between two expected $L_2-$discrepancies, thus

\begin{equation}\label{eq438}
\begin{aligned}
    &\mathbb{E}L_2^2(P_{\Omega^{*}_{\sim}})-\mathbb{E}L_2^2(P_{\Omega^{*}_{|}})
    \\&=\frac{1}{N^2}\sum_{i=1}^{2}\int_{I_{d,m}^{0}}(\tilde{q}_i(x)-\bar{q}_i(x))dx+\frac{1}{N^2}\sum_{i=1}^{2}\int_{I_{d,m}^{0}}(\bar{q}^2_i(x)-\tilde{q}^2_i(x))dx\\&=\frac{1}{N}\sum_{i=1}^{2}\int_{I_{d,m}^{0}}(\lambda(\tilde{\Omega}_i\cap[0,x])-\lambda(\bar{\Omega}_i\cap[0,x]))dx\\&+\sum_{i=1}^{2}\int_{I_{d,m}^{0}}(\lambda^2(\bar{\Omega}_i\cap[0,x])-\lambda^2(\tilde{\Omega}_i\cap[0,x]))dx\\&=\frac{1}{N^3}\sum_{i=1}^{2}\int_{I'_d}(\lambda(\Omega''_{i,\sim}\cap[0,x])-\lambda(\Omega''_{i,|}\cap[0,x]))dx\\&+\frac{1}{N^3}\sum_{i=1}^{2}\int_{I'_d}(\lambda^2(\Omega''_{i,|}\cap[0,x])-\lambda^2(\Omega''_{i,\sim}\cap[0,x]))dx.
\end{aligned}
\end{equation}

Furthermore, employing \eqref{el22pom} again, we have

\begin{equation}\label{eq439}
    \begin{aligned}
    &\mathbb{E}(L_2^2(P_{\Omega''_{\sim}}))-\mathbb{E}(L_2^2(P_{\Omega''_{|}}))\\&=\frac{1}{8}\sum_{i=1}^{2}\int_{I'_d}\mathbf{q}'_{i,\sim}(\mathbf{x})(1-\mathbf{q}'_{i,\sim}(\mathbf{x}))dx-\frac{1}{8}\sum_{i=1}^{2}\int_{I'_d}\mathbf{q}'_{i,|}(\mathbf{x})(1-\mathbf{q}'_{i,|}(\mathbf{x}))dx\\&=\frac{1}{8}\sum_{i=1}^{2}\int_{I'_d}(\lambda(\Omega''_{i,\sim}\cap[0,x])-\lambda(\Omega''_{i,|}\cap[0,x]))dx\\&+\frac{1}{8}\sum_{i=1}^{2}\int_{I'_d}(\lambda^2(\Omega''_{i,|}\cap[0,x])-\lambda^2(\Omega''_{i,\sim}\cap[0,x]))dx.
    \end{aligned}
\end{equation}

Combining with \eqref{eq438} and \eqref{eq439}, we obtain

\begin{equation}\label{eq431}
\mathbb{E}L_2^2(P_{\Omega^{*}_{\sim}})-\mathbb{E}L_2^2(P_{\Omega^{*}_{|}})=\frac{1}{N^3}\cdot [8\mathbb{E}(L_2^2(P_{\Omega''_{\sim}}))-8\mathbb{E}(L_2^2(P_{\Omega''_{|}}))].
\end{equation}

Combining with \eqref{b1simb2sim}, \eqref{b1simb2sim111} and \eqref{8l22omesim}, the proof is completed.

\subsection{Proof of Theorem \ref{uniformnoise0}}
Following the proof process of Lemma \ref{lm43}, we obtain Theorem \ref{uniformnoise0}.

\subsection{Proof of Theorem \ref{expstarnp}}
We only consider the case $arctan\frac{1}{2}\leq\theta\leq\frac{\pi}{2}$, the calculation of case $0\leq\theta<arctan\frac{1}{2}$ is similar to it. First, we have
$$
    P_2(\theta)=-\frac{1}{24tan\theta}+\frac{1}{120tan^2\theta}+\frac{1}{1440tan^3\theta}.
$$ 

Then from Lemma \ref{lm43}, we obtain

\begin{equation}\label{diffepn}
    \mathbb{E}L_2^2(P_{\Omega^{*}_{\sim}})-\mathbb{E}L_2^2(P_{\Omega^{*}_{|}})=\frac{1}{N^3}\cdot\frac{1}{3^{d-2}}\cdot P_{2}(\theta),
\end{equation}
where 
$$
    P_{\Omega^{*}_{\sim}}=\{U_1,U_2,\ldots,U_N\},
$$
and

$$
    P_{\Omega^{*}_{|}}=\{W_1,W_2,\ldots,W_N\},
$$
denote stratified samples under different partition models $\Omega^{*}_{\sim}$ and $\Omega^{*}_{|}$ respectively.

Now, for arbitrary test set $R=[0,x)\subset [0,1]^d$, we consider the following discrepancy function,

\begin{equation}\label{dispfunc1}
\Delta_{\mathscr{P}}(x)=\frac{1}{N}
\sum_{n=1}^{N}\mathbf{1}_{R}(W_n)-\lambda(R).
\end{equation}

For an equivolume partition $\Omega=\{\Omega_1,\Omega_2,\ldots,\Omega_N\}$, we divide the test set $R$ into two parts, one is the disjoint union of $\Omega_{i}$ entirely contained by $R$ and another is the union of remaining pieces which are the intersections of some $\Omega_{j}$ and $R$, i.e.,

\begin{equation}\label{Rtp1}
R=\bigcup_{i\in I_0}\Omega_{i}\cup\bigcup_{j\in J_0}(\Omega_{j}\cap R),
\end{equation}
where $I_0,J_0$ are two index-sets.

Let $$T=\bigcup_{j\in J_0}(\Omega_{j}\cap R),$$ then from \eqref{dispfunc1}, we have

\begin{equation}\label{dfots}
\Delta_{\mathscr{P}}(x)=\frac{1}{N}
\sum_{n=1}^{N}\mathbf{1}_{R}(W_n)-\lambda(R)=\frac{1}{N}
\sum_{n=1}^{N}\mathbf{1}_{T}(W_n)-\lambda(T),
\end{equation}
where \eqref{dfots} is based on the fact discrepancy function equals $0$ on $\bigcup_{i\in I_0}\Omega_{i}$.

According to the definition of $L_2-$discrepancy and \eqref{dfots}, it follows that

\begin{equation}\label{exl2dpb1_1}
\mathbb{E}(L_2^2(D_N,P_{\Omega^{*}_{|}}))=\mathbb{E}
(\int_{[0,1]^d}|\frac{1}{N}
\sum_{n=1}^{N}\mathbf{1}_{T}(W_n)-\lambda(T)|^2dx).
\end{equation}

Consider the whole sum in \eqref{exl2dpb1_1} as a random variable which is defined on a region $P_{\Omega}$. Besides we set the probability measure be $w$, then we have

\begin{equation}\label{el22dnz}
\begin{aligned}
\mathbb{E}(L_2^2(D_N,P_{\Omega^{*}_{|}}))&=\int_{P_{\Omega}}
\int_{[0,1]^d}|\frac{1}{N}
\sum_{n=1}^{N}\mathbf{1}_{T}(W_n)-\lambda(T)|^2dxdw\\&
=\int_{[0,1]^d}\int_{P_{\Omega}}|\frac{1}{N}
\sum_{n=1}^{N}\mathbf{1}_{T}(W_n)-\lambda(T)|^2dwdx.
\end{aligned}
\end{equation}

It can easily be checked that,

$$
\mathbb{E}(\frac{1}{N}
\sum_{n=1}^{N}\mathbf{1}_{T}(W_n))=\int_{P_{\Omega}}\frac{1}{N}
\sum_{n=1}^{N}\mathbf{1}_{T}(W_n)dw=\lambda(T).
$$

Hence,

\begin{equation}\label{expvar1}
\int_{P_{\Omega}}|\frac{1}{N}
\sum_{n=1}^{N}\mathbf{1}_{T}(W_n)-\lambda(T)|^2dw=\text{Var}(\frac{1}{N}\sum_{n=1}^{N}\mathbf{1}_{T}(W_n)).
\end{equation}

Let $\sigma_{j}^{2}=\text{Var}(\mathbf{1}_{T}(W_j)),$ then we have

\begin{equation}
    \Sigma^{2}=\text{Var}(\sum_{n=1}^{N}\mathbf{1}_{T}(W_n))=\sum_{n=1}^{N}\text{Var}(\mathbf{1}_{T}(W_n))=\sum\limits_{j\in J_0}\sigma_{j}^{2}.
\end{equation}

Hence, from \eqref{CNbd0}, we get

\begin{equation}\label{sgm2}
\Sigma^{2}\leq d\cdot N^{1-\frac{1}{d}}.
\end{equation}

Therefore,

\begin{equation}\label{el22pgp}
\begin{aligned}
    \mathbb{E}L_2^2(P_{\Omega^{*}_{|}})&=\int_{[0,1]^{d}}\text{Var}(\frac{1}{N}\sum_{n=1}^{N}\mathbf{1}_{T}(W_n))dx\\&\leq \frac{d}{N^{1+\frac{1}{d}}}.
\end{aligned}
\end{equation}

Combining with \eqref{diffepn} and \eqref{el22pgp}, the desired result is proved.

\section{Conclusion}\label{conclu}

We study  expected $L_2-$discrepancy under a class of new convex equivolume partitions. First, the expected $L_2-$discrepancy under two partition models are compared. Second, the explicit expected $L_2-$discrepancy upper bounds under the new partition models are obtained. So the optimal partition model that minimizes expected $L_2-$discrepancy is found and an optimal expected $L_2-$discrepancy upper bound is given explicitly under a class of new convex equivolume partitions. In future, star discrepancy will be studied under a class of convex equal volume partitions, which will have more corresponding applications.

\end{document}